\documentclass[12pt]{article}
\usepackage{amsmath,amssymb,amsthm}
\numberwithin{equation}{section}
\usepackage{units}
\usepackage{color}
\usepackage[T1]{fontenc}
\usepackage[utf8]{inputenc}
\usepackage{authblk}
\usepackage{bm}
\usepackage{graphicx}

\usepackage{datetime}

\usepackage[colorlinks=true,
linkcolor=webgreen,
filecolor=webbrown,
citecolor=webgreen]{hyperref}

\definecolor{webgreen}{rgb}{0,.5,0}
\definecolor{webbrown}{rgb}{.6,0,0}

\usepackage[toc,page]{appendix}

\newcommand{\Z}{{\mathbb Z}}

\newtheorem{thm}{Theorem}
\newtheorem{theorem}[thm]{Theorem}
\newtheorem{lemma}{Lemma}

\newtheorem{corollary}[thm]{Corollary}

\title{An Identity for Second Order Sequences Obeying the Same Recurrence Relation}

\author[]{Kunle Adegoke \\\href{mailto:adegoke00@gmail.com}{\tt adegoke00@gmail.com}}

\affil{Department of Physics and Engineering Physics, \mbox{Obafemi Awolowo University}, 220005 Ile-Ife, Nigeria}

\begin{document}

\date{}

\maketitle

\begin{abstract}
\noindent We derive an identity connecting any two second-order linear recurrence sequences having the same recurrence relation but whose initial terms may be different. Binomial and ordinary summation identities arising from the identity are developed. Illustrative examples are drawn from Fibonacci, Fibonacci-Lucas, Pell, Pell-Lucas, Jacobsthal and Jacobsthal-Lucas sequences and their generalizations. Our new results subsume previously known identities.

\end{abstract}
\section{Introduction}
This paper is concerned with establishing an identity connecting any two second-order homogeneous linear recurrence sequences, $(G_n)_{n\in\Z}$ and $(H_n)_{n\in\Z}$, having the same recurrence relation but whose initial terms may be different. Thus, for $n\ge2$ and with $p$ and $q$ arbitrary fixed non-zero complex constants, we define 
\begin{equation}
G_n=pG_{n-1}+qG_{n-2}\,,
\end{equation}
where the initial terms $G_0$ and $G_1$ are given arbitrary integers, not both zero; and 
\begin{equation}
H_n=pH_{n-1}+qH_{n-2}\,,
\end{equation}
with initial terms $H_0$ and $H_1$ given arbitrary integers, not both zero.

\medskip

Extension of the definition of $G_n$ and $H_n$ to negative subscripts is provided by writing the recurrence relation as
\begin{equation}
G_{-n}=(G_{-n+2}-pG_{-n+1})/q
\end{equation}
and
\begin{equation}
H_{-n}=(H_{-n+2}-pH_{-n+1})/q\,.
\end{equation}
In section \ref{sec.main}, we will derive an identity connecting $(G_n)$ and $(H_n)$, for arbitrary integers. We will illustrate the results by deriving identities for six well-known second-order sequences, namely, Fibonacci, Lucas, Pell, Pell-Lucas, Jacobsthal and Jacobsthal-Lucas sequences.
\subsection{Notation and definitions}
The Fibonacci numbers, $F_n$, and the Lucas numbers, $L_n$, are defined, for \mbox{$n\in\Z$}, as usual, through the recurrence relations \mbox{$F_n=F_{n-1}+F_{n-2}$ ($n\ge 2$)}, \mbox{$F_0=0$, $F_1=1$} and \mbox{$L_n=L_{n-1}+L_{n-2}$} ($n\ge 2$), $L_0=2$, $L_1=1$, with $F_{-n}=(-1)^{n-1}F_n$ and $L_{-n}=(-1)^nL_n$. Exhaustive discussion of the properties of Fibonacci and Lucas numbers can be found in Vajda~\cite{vajda} and in Koshy~\cite{koshy}. Generalized Fibonacci numbers having the same recurrence as the Fibonacci and Lucas numbers but with arbitrary initial values will be denoted $\mathcal{F}_n$.

\medskip

The Jacobsthal numbers, $J_n$, and the Jacobsthal-Lucas numbers, $j_n$, are defined, for \mbox{$n\in\Z$}, through the recurrence relations \mbox{$J_n=J_{n-1}+2J_{n-2}$} ($n\ge 2$), \mbox{$J_0=0$, $J_1=1$} and \mbox{$j_n=j_{n-1}+2j_{n-2}$} ($n\ge 2$), $j_0=2$, $j_1=1$, with $J_{-n}=(-1)^{n-1}2^{-n}J_n$ and $j_{-n}=(-1)^n2^{-n}j_n$. Horadam~\cite{horadam96} and Aydin~\cite{aydin} are good reference materials on the Jacobsthal and associated sequences.  Generalized Jacobsthal numbers having the same recurrence as the Jacobsthal and Jacbsthal-Lucas numbers but with arbitrary initial values will be denoted $\mathcal{J}_n$.

\medskip

The Pell numbers, $P_n$, and Pell-Lucas numbers, $Q_n$, are defined, for \mbox{$n\in\Z$}, through the recurrence relations \mbox{$P_n=2P_{n-1}+P_{n-2}$} ($n\ge 2$), \mbox{$P_0=0$, $P_1=1$} and \mbox{$Q_n=2Q_{n-1}+Q_{n-2}$} ($n\ge 2$), $Q_0=2$, $Q_1=2$, with $P_{-n}=(-1)^{n-1}P_n$ and $Q_{-n}=(-1)^nQ_n$. Koshy~\cite{koshy14}, Horadam~\cite{horadam71} and Patel and Shrivastava~\cite{patel13} are useful source materials on Pell and Pell-Lucas numbers.  Generalized Pell numbers having the same recurrence as the Pell and Pell-Lucas numbers but with arbitrary initial values will be denoted $\mathcal{P}_n$.

\medskip

For reference, the first few values of the six sequences are given below:
\begin{center}
\begin{tabular}{cllllllllllllll}
$n$: & \multicolumn{1}{c}{$-5$} & \multicolumn{1}{c}{$ -4$} & \multicolumn{1}{c}{$ -3$} & \multicolumn{1}{c}{$ -2$} & \multicolumn{1}{c}{$ -1$} & \multicolumn{1}{c}{$ 0$} & \multicolumn{1}{c}{$ 1$} & \multicolumn{1}{c}{$ 2$} & \multicolumn{1}{c}{$ 3$} & \multicolumn{1}{c}{$ 4$} & \multicolumn{1}{c}{$ 5$} & \multicolumn{1}{c}{$ 6$} & \multicolumn{1}{c}{$ 7$} & \multicolumn{1}{c}{$ 8$} \\ 
\hline
$F_n$: & \multicolumn{1}{c}{$5$} & \multicolumn{1}{c}{$ -3$} & \multicolumn{1}{c}{$ 2$} & \multicolumn{1}{c}{$ -1$} & \multicolumn{1}{c}{$ 1$} & \multicolumn{1}{c}{$ 0$} & \multicolumn{1}{c}{$ 1$} & \multicolumn{1}{c}{$ 1$} & \multicolumn{1}{c}{$ 2$} & \multicolumn{1}{c}{$ 3$} & \multicolumn{1}{c}{$ 5$} & \multicolumn{1}{c}{$ 8$} & \multicolumn{1}{c}{$ 13$} & \multicolumn{1}{c}{$ 21$} \\ 
$L_n$: & \multicolumn{1}{c}{$-11$} & \multicolumn{1}{c}{$ 7$} & \multicolumn{1}{c}{$ -4$} & \multicolumn{1}{c}{$ 3$} & \multicolumn{1}{c}{$ -1$} & \multicolumn{1}{c}{$ 2$} & \multicolumn{1}{c}{$ 1$} & \multicolumn{1}{c}{$ 3$} & \multicolumn{1}{c}{$ 4$} & \multicolumn{1}{c}{$ 7$} & \multicolumn{1}{c}{$ 11$} & \multicolumn{1}{c}{$ 18$} & \multicolumn{1}{c}{$ 29$} & \multicolumn{1}{c}{$ 47$} \\ 
$P_n$: & \multicolumn{1}{c}{$29$} & \multicolumn{1}{c}{$ -12$} & \multicolumn{1}{c}{$ 5$} & \multicolumn{1}{c}{$ -2$} & \multicolumn{1}{c}{$ 1$} & \multicolumn{1}{c}{$ 0$} & \multicolumn{1}{c}{$ 1$} & \multicolumn{1}{c}{$ 2$} & \multicolumn{1}{c}{$ 5$} & \multicolumn{1}{c}{$ 12$} & \multicolumn{1}{c}{$ 29$} & \multicolumn{1}{c}{$ 70$} & \multicolumn{1}{c}{$ 169$} & \multicolumn{1}{c}{$ 408$} \\ 
$Q_n$: & \multicolumn{1}{c}{$-82$} & \multicolumn{1}{c}{$ 34$} & \multicolumn{1}{c}{$ -14$} & \multicolumn{1}{c}{$ 6$} & \multicolumn{1}{c}{$ -2$} & \multicolumn{1}{c}{$ 2$} & \multicolumn{1}{c}{$ 2$} & \multicolumn{1}{c}{$ 6$} & \multicolumn{1}{c}{$ 14$} & \multicolumn{1}{c}{$ 34$} & \multicolumn{1}{c}{$ 82$} & \multicolumn{1}{c}{$ 198$} & \multicolumn{1}{c}{$ 478$} & \multicolumn{1}{c}{$ 1154$} \\ 
$J_n$: & \multicolumn{1}{c}{$11/32$} & \multicolumn{1}{c}{$ -5/16$} & \multicolumn{1}{c}{$ 3/8$} & \multicolumn{1}{c}{$ -1/4$} & \multicolumn{1}{c}{$ 1/2$} & \multicolumn{1}{c}{$ 0$} & \multicolumn{1}{c}{$ 1$} & \multicolumn{1}{c}{$ 1$} & \multicolumn{1}{c}{$ 3$} & \multicolumn{1}{c}{$ 5$} & \multicolumn{1}{c}{$ 11$} & \multicolumn{1}{c}{$ 21$} & \multicolumn{1}{c}{$ 43$} & \multicolumn{1}{c}{$ 85$} \\ 
$j_n$: & \multicolumn{1}{c}{$-31/32$} & \multicolumn{1}{c}{$ 17/16$} & \multicolumn{1}{c}{$ -7/8$} & \multicolumn{1}{c}{$ 5/4$} & \multicolumn{1}{c}{$ -1/2$} & \multicolumn{1}{c}{$ 2$} & \multicolumn{1}{c}{$ 1$} & \multicolumn{1}{c}{$ 5$} & \multicolumn{1}{c}{$ 7$} & \multicolumn{1}{c}{$ 17$} & \multicolumn{1}{c}{$ 31$} & \multicolumn{1}{c}{$ 65$} & \multicolumn{1}{c}{$ 127$} & \multicolumn{1}{c}{$ 257$} \\ 
\end{tabular}
\end{center}

\section{Main results}\label{sec.main}
\subsection{Recurrence relations and related identities}
\begin{theorem}\label{thm.ash3b4f}
Let $(G_n)_{n\in\Z}$ and $(H_n)_{n\in\Z}$ be any two linear homogeneous recurrence sequences having the same recurrence relation. Then, the following identity holds for arbitrary integers $n$, $m$, $a$, $b$, $c$ and $d$:
\[
\begin{split}
&(G_{d-b}G_{c-a}-G_{d-a}G_{c-b})H_{n+m}\\
&\qquad\qquad= (G_{d-b}G_{m-a}-G_{d-a}G_{m-b})H_{n+c}\\
&\quad\qquad\qquad+(G_{c-a}G_{m-b}-G_{c-b}G_{m-a})H_{n+d}\,.
\end{split}
\]
\end{theorem}
In particular, we have
\begin{equation}\label{eq.nveox2i}
\begin{split}
&(G_{d-b}G_{c-a}-G_{d-a}G_{c-b})G_{n+m}\\
&\qquad\qquad= (G_{d-b}G_{m-a}-G_{d-a}G_{m-b})G_{n+c}\\
&\quad\qquad\qquad+(G_{c-a}G_{m-b}-G_{c-b}G_{m-a})G_{n+d}\,,
\end{split}
\end{equation}
for any linear, second-order recurrence sequence, $(G_n)_{n\in\Z}$.
\begin{proof}
Since both sequences $(G_n)$ and $(H_n)$ have the same recurrence relation we choose a basis set in one and express the other in this basis. We write
\begin{equation}\label{eq.nbl33qg}
H_{n+m}=\lambda_1G_{m-a}+\lambda_2G_{m-b}\,,
\end{equation}
where $a$, $b$, $n$ and $m$ are arbitrary integers and the coefficients $\lambda_1$ and $\lambda_2$ are to be determined. Setting $m=c$ and $m=d$, in turn, produces two simultaneous equations:
\[
H_{n + c}  = \lambda_1 G_{c - a}  + \lambda_2 G_{c - b},\quad H_{n + d}  = \lambda_1 G_{d - a}  + \lambda_2 G_{d - b}\,.
\]
The identity of Theorem \ref{thm.ash3b4f} is established by solving these equations for $\lambda_1$ and $\lambda_2$ and substituting the solutions into identity \eqref{eq.nbl33qg}.
\end{proof}
\begin{corollary}\label{cor.g869wix}
The following identity holds for integers $a$, $b$, $n$ and $m$:
\[
\begin{split}
&(G_{a-b}G_{b-a}-G_{0}^2\,)H_{n+m}\\ 
&\qquad= (G_{b-a}G_{m-b}-G_{0}G_{m-a})H_{n+a}\\
&\qquad\quad+(G_{a-b}G_{m-a}-G_{0}G_{m-b})H_{n+b}\,.
\end{split}
\]
\end{corollary}
In particular,
\begin{equation}
\begin{split}
&(G_{a-b}G_{b-a}-G_{0}^2\,)G_{n+m}\\ 
&\qquad= (G_{b-a}G_{m-b}-G_{0}G_{m-a})G_{n+a}\\
&\qquad\quad+(G_{a-b}G_{m-a}-G_{0}G_{m-b})G_{n+b}\,.
\end{split}
\end{equation}
\subsection{Summation identities}
\subsubsection{Summation identities not involving binomial coefficients}
\begin{lemma}[{\cite[Lemma 1]{adegoke18}}]\label{lem.u4bqbkc}
Let $\{X_n\}$ and $\{Y_n\}$ be any two sequences such that $X_n$ and $Y_n$, $n\in\Z$, are connected by a three-term recurrence relation $X_n=f_1X_{n-a}+f_2Y_{n-b}$, where $f_1$ and $f_2$ are arbitrary non-vanishing complex functions, not dependent on $r$, and $a$ and $b$ are integers. Then,
\[
f_2\sum_{j = 0}^k {\frac{{Y_{n - ka  - b  + a j} }}{{f_1{}^j }}}  = \frac{{X_n }}{{f_1{}^k }} - f_1X_{n - (k + 1)a }\,, 
\]
for $k$ a non-negative integer.
\end{lemma}

\begin{lemma}[{\cite[Lemma 2]{adegoke18}}]\label{lem.s9jfs7n}
Let $\{X_n\}$ be any arbitrary sequence, where $X_n$, $n\in\Z$, satisfies a three-term recurrence relation $X_n=f_1X_{n-a}+f_2X_{n-b}$, where $f_1$ and $f_2$ are arbitrary non-vanishing complex functions, not dependent on $r$, and $a$ and $b$ are integers. Then, the following identities hold for integer $k$:
\begin{equation}\label{eq.mxyb9zk}
f_2\sum_{j = 0}^k {\frac{{X_{n - ka  - b  + a j} }}{{f_1^j }}}  = \frac{{X_n }}{{f_1^k }} - f_1X_{n - (k + 1)a }\,,
\end{equation}
\begin{equation}\label{eq.cgldajj}
f_1\sum_{j = 0}^k {\frac{{X_{n - kb  - a  + b j} }}{{f_2^j }}}  = \frac{{X_n }}{{f_2^k }} - f_2X_{n - (k + 1)b } 
\end{equation}
and
\begin{equation}\label{eq.n2n4ec3}
\sum_{j = 0}^k { \frac{X_{n - (a - b)k + b + (a - b)j}}{(-f_1/f_2)^j} }  = \frac{f_2 X_n}{(-f_1/f_2)^k}  + f_1 X_{n - (k + 1)(a - b)}\,.
\end{equation}

\end{lemma}
In order to state the next two theorems in a compact form, we introduce the following notation:
\begin{equation}
f_G(u,v;s,t)=G_{u-s}G_{v-t}-G_{u-t}G_{v-s}\,,
\end{equation}
with the following symmetry properties:
\begin{equation}
f_G(u,v;t,s)=-f_G(u,v;s,t)\,,\quad f_G(v,u;s,t)=-f_G(u,v;s,t)\,,
\end{equation}
\begin{equation}
f_G(v,u;t,s)=f_G(u,v;s,t)\,,
\end{equation}
and
\begin{equation}
f_G(u,u;s,t)=0\,,\quad f_G(u,v;s,s)=0\,.
\end{equation}
In this notation, the identity of Theorem \ref{thm.ash3b4f} becomes
\begin{equation}\label{eq.kjhpj9e}
f_G(d,c;b,a)H_{n + m}  = f_G(d,m;b,a)H_{n + c}  + f_G(c,m;a,b)H_{n + d}\,.
\end{equation}
The results in the next theorem follow from direct substitutions from identity \eqref{eq.kjhpj9e} into Lemma \ref{lem.s9jfs7n}.
\begin{theorem}
The following identities hold for arbitrary integers $a$, $b$, $c$, $d$ and $m$ for which $f_G(d,m;b,a)\ne 0$ and $f_G(c,m;b,a)\ne 0$:
\begin{equation}
\begin{split}
&{f_G(c,m;a,b)}\sum\limits_{j = 0}^k {\left( {\frac{{f_G(d,c;b,a)}}{{{f_G(d,m;b,a)}}}} \right)^j H_{n - (m-c) k - (m-d)  + (m-c) j} }\\
&\qquad\qquad\qquad= \frac{{{f_G(d,c;b,a)}^{k + 1} }}{{{f_G(d,m;b,a)}^k }}H_n  - {f_G(d,m;b,a)}H_{n - (m-c) (k + 1)} \,,
\end{split}
\end{equation}
\begin{equation}
\begin{split}
&{f_G(d,m;b,a)}\sum\limits_{j = 0}^k {\left( {\frac{{f_G(d,c;b,a)}}{{{f_G(c,m;a,b)}}}} \right)^j H_{n - (m-d) k - (m-c)  + (m-d) j} }\\
&\qquad\qquad\qquad= \frac{{{f_G(d,c;b,a)}^{k + 1} }}{{{f_G(c,m;a,b)}^k }}H_n  - {f_G(c,m;a,b)}H_{n - (m-d) (k + 1)} 
\end{split}
\end{equation}
and
\begin{equation}
\begin{split}
&{f_G(d,c;b,a)}\sum\limits_{j = 0}^k {\left( { - \frac{{{f_G(d,m;b,a)}}}{{{f_G(c,m;a,b)}}}} \right)^j H_{n - (c-d) k + (m-c)  + (c-d) j} }\\
&\qquad\qquad\qquad = ( - 1)^k \frac{{{f_G(d,m;b,a)}^{k + 1} }}{{{f_G(c,m;a,b)}^k }}H_n  + {f_G(c,m;a,b)}H_{n - (c-d) (k + 1)}\,.
\end{split}
\end{equation}
\end{theorem}
\subsubsection{Binomial summation identities}
\begin{lemma}[{\cite[Lemma 3]{adegoke18}}]\label{lem.binomial}
Let $\{X_n\}$ be any arbitrary sequence. Let $X_n$, $n\in\Z$, satisfy a three-term recurrence relation $X_n=f_1X_{n-a}+f_2X_{n-b}$, where $f_1$ and $f_2$ are non-vanishing complex functions, not dependent on $n$, and $a$ and $b$ are integers. Then,
\begin{equation}\label{eq.fe496kc}
\sum_{j = 0}^k {\binom kj\left( {\frac{f_1}{f_2}} \right)^j X_{n - bk  + (b  - a )j} }  = \frac{{X_n }}{{f_2^k }}\,,
\end{equation}
\begin{equation}\label{eq.j7k6a8g}
\sum_{j = 0}^k {\binom kj\frac{{X_{n + (a - b)k + bj} }}{{( - f_2 )^j }}}  = \left( { - \frac{{f_1 }}{{f_2 }}} \right)^k X_n
\end{equation}
and
\begin{equation}\label{eq.fnwrzi3}
\sum_{j = 0}^k {\binom kj\frac{{X_{n + (b - a)k + a j} }}{{( - f_1 )^j }}}  = \left( { - \frac{f_2}{f_1}} \right)^k X_n\,,
\end{equation}
for $k$ a non-negative integer.

\end{lemma}
Substituting from identity \eqref{eq.kjhpj9e} into Lemma \ref{lem.binomial}, we have the results stated in the next theorem.
\begin{theorem}
The following identities hold for positive integer $k$ and arbitrary integers $a$, $b$, $c$, $d$ and $m$ for which $f_G(d,m;b,a)\ne 0$ and $f_G(c,m;b,a)\ne 0$:
\begin{equation}
\sum\limits_{j = 0}^k {\binom kj\left( {\frac{{{f_G(d,m;b,a)}}}{{{f_G(c,m;a,b)}}}} \right)^j H_{n - (m-d) k + (c-d) j} }  = \left( {\frac{f_G(d,c;b,a)}{{{f_G(c,m;a,b)}}}} \right)^k H_n\,,
\end{equation}
\begin{equation}
\sum\limits_{j = 0}^k {\binom kj\left( { - \frac{f_G(d,c;b,a)}{{{f_G(c,m;a,b)}}}} \right)^j H_{n - (c-d) k + (m-d) j} }  = \left( { - \frac{{{f_G(d,m;b,a)}}}{{{f_G(c,m;a,b)}}}} \right)^k H_n
\end{equation}
and
\begin{equation}
\sum\limits_{j = 0}^k {\binom kj\left( { - \frac{f_G(d,c;b,a)}{{{f_G(d,m;b,a)}}}} \right)^j H_{n + (c-d) k + (m-c) j} }  = \left( { - \frac{{{f_G(c,m;a,b)}}}{{{f_G(d,m;b,a)}}}} \right)^k H_n\,.
\end{equation}

\end{theorem}
\section{Explicit examples from well-known second order sequences}
\subsection{Identities involving Fibonacci, Lucas and generalized Fibonacci numbers}\label{sec.fibonacci}
\subsubsection{Recurrence relations and related identities}
In Corollary \ref{cor.g869wix}, let $(G_n)\equiv(F_n)$ be the Fibonacci sequence and let $(H_n)\equiv(\mathcal{F}_n)$ be the generalized Fibonacci sequence. Then, the identity of Corollary \ref{cor.g869wix} reduces to
\begin{equation}\label{eq.aa1gc12}
F_{a-b}\mathcal{F}_{n+m}=F_{m-b}\mathcal{F}_{n+a}-(-1)^{a-b}F_{m-a}\mathcal{F}_{n+b}\,.
\end{equation}
The presumably new identity \eqref{eq.aa1gc12} subsumes most known three-term recurrence relations involving Fibonacci numbers, Lucas numbers and the generalized Fibonacci numbers. We will give a couple of examples to illustrate this point.

\medskip

Incidentally, identity \eqref{eq.aa1gc12} can also be written as
\begin{equation}\label{eq.h84fbk1}
F_{a-b}\mathcal{F}_{n+m}=\mathcal{F}_{m-b}F_{n+a}-(-1)^{a-b}\mathcal{F}_{m-a}F_{n+b}\,.
\end{equation}
Setting $a=0$, $b=m-n$ in identity \eqref{eq.aa1gc12} gives
\begin{equation}
F_{n-m}\mathcal{F}_{n+m}=F_n\mathcal{F}_n-(-1)^{n-m}F_m\mathcal{F}_m\,,
\end{equation}
which is a generalization of Catalan's identity:
\begin{equation}
F_{n-m}F_{n+m}=F_n^2+(-1)^{n+m+1}F_m^2\,.
\end{equation}
Setting $b=-a$ in identity \eqref{eq.aa1gc12} gives
\begin{equation}
F_{2a}\mathcal{F}_{n+m}=F_{m+a}\mathcal{F}_{n+a}-F_{m-a}\mathcal{F}_{n-a}\,,
\end{equation}
with the special case
\begin{equation}
\mathcal{F}_{n+m}=F_{m+1}\mathcal{F}_{n+1}-F_{m-1}\mathcal{F}_{n-1}\,,
\end{equation}
which is a generalization of the following known identity (Halton \cite[Identity (63)]{halton65}):
\begin{equation}
F_{n+m}=F_{m+1}F_{n+1}-F_{m-1}F_{n-1}\,.
\end{equation}
Setting $b=2k$, $a=1$ and $b=2k$, $a=0$, in turn, in identity \eqref{eq.aa1gc12} produces
\begin{equation}
F_{2k - 1} \mathcal{F}_{n + m}  = F_{m - 2k} \mathcal{F}_{n + 1}  + F_{m - 1} \mathcal{F}_{n + 2k}
\end{equation}
and
\begin{equation}\label{eq.m2jtlr3}
F_{2k} \mathcal{F}_{n + m}  = F_m \mathcal{F}_{n + 2k}  - F_{m - 2k} \mathcal{F}_n\,.
\end{equation}
Identity \eqref{eq.m2jtlr3} is a generalization of the known identity (Vajda \cite[Formula (8)]{vajda}):
\begin{equation}
\mathcal{F}_{n + m}  = F_{m - 1} \mathcal{F}_{n}  + F_{m} \mathcal{F}_{n + 1}\,.
\end{equation}
Setting $a=n$ and $b=-m$ in \eqref{eq.aa1gc12} produces
\begin{equation}
F_{2m} \mathcal{F}_{2n}  = F_{n + m} \mathcal{F}_{n + m}  - F_{n - m} \mathcal{F}_{n - m}\,.
\end{equation}
\subsubsection{Summation identities}
Summation identities given in Theorems \ref{thm.ik24j18} and \ref{thm.tutbc4x} are derived by making appropriate substitutions from identity \eqref{eq.aa1gc12} into Lemmata \ref{lem.s9jfs7n} and \ref{lem.binomial}. 
\begin{theorem}\label{thm.ik24j18}
The following identities hold for arbitrary integers $a$, $b$, $n$, $m$ and $k$:
\begin{equation}
\begin{split}
&( - 1)^{a + b + 1} F_{m - a} \sum\limits_{j = 0}^k {F_{m - b} ^{k - j} F_{a - b} ^j \mathcal{F}_{n - (m - a)k - (m - b) + (m - a)j} }\\
&\qquad\qquad\qquad\qquad= F_{a - b} ^{k + 1} \mathcal{F}_n  - F_{m - b} ^{k + 1} \mathcal{F}_{n - (m - a)(k + 1)}\,,
\end{split}
\end{equation}
\begin{equation}
\begin{split}
&F_{m - b} \sum\limits_{j = 0}^k {( - 1)^{(a + b + 1)(k - j)} F_{m - a} ^{k - j} F_{a - b} ^j \mathcal{F}_{n - (m - b)k - (m - a) + (m - b)j} }\\
&\qquad\qquad\qquad= F_{a - b} ^{k + 1} \mathcal{F}_n  - ( - 1)^{(a + b + 1)(k + 1)} F_{m - a} ^{k + 1} \mathcal{F}_{n - (m - b)(k + 1)}
\end{split}
\end{equation}
and
\begin{equation}
\begin{split}
&F_{a - b} \sum\limits_{j = 0}^k {( - 1)^{(a + b)j} F_{m - a} ^{k - j} F_{m - b} ^j \mathcal{F}_{n - (a - b)k + (m - a) + (a - b)j} }\\
&\qquad\qquad= ( - 1)^{(a + b)k} F_{m - b} ^{k + 1} \mathcal{F}_n  + ( - 1)^{a + b + 1} F_{m - a} ^{k + 1} \mathcal{F}_{n - (a - b)(k + 1)}\,.
\end{split}
\end{equation}

\end{theorem}
In particular, we have
\begin{equation}
\begin{split}
&( - 1)^{a + b + 1} F_{m - a} \sum\limits_{j = 0}^k {F_{m - b} ^{k - j} F_{a - b} ^j F_{n - (m - a)k - (m - b) + (m - a)j} }\\
&\qquad\qquad\qquad\qquad= F_{a - b} ^{k + 1} F_n  - F_{m - b} ^{k + 1} F_{n - (m - a)(k + 1)}\,,
\end{split}
\end{equation}
\begin{equation}
\begin{split}
&F_{m - b} \sum\limits_{j = 0}^k {( - 1)^{(a + b + 1)(k - j)} F_{m - a} ^{k - j} F_{a - b} ^j F_{n - (m - b)k - (m - a) + (m - b)j} }\\
&\qquad\qquad\qquad= F_{a - b} ^{k + 1} F_n  - ( - 1)^{(a + b + 1)(k + 1)} F_{m - a} ^{k + 1} F_{n - (m - b)(k + 1)}
\end{split}
\end{equation}
and
\begin{equation}
\begin{split}
&F_{a - b} \sum\limits_{j = 0}^k {( - 1)^{(a + b)j} F_{m - a} ^{k - j} F_{m - b} ^j F_{n - (a - b)k + (m - a) + (a - b)j} }\\
&\qquad\qquad= ( - 1)^{(a + b)k} F_{m - b} ^{k + 1} F_n  + ( - 1)^{a + b + 1} F_{m - a} ^{k + 1} F_{n - (a - b)(k + 1)}\,;
\end{split}
\end{equation}
and the corresponding results involving Lucas numbers:
\begin{equation}
\begin{split}
&( - 1)^{a + b + 1} F_{m - a} \sum\limits_{j = 0}^k {F_{m - b} ^{k - j} F_{a - b} ^j L_{n - (m - a)k - (m - b) + (m - a)j} }\\
&\qquad\qquad\qquad\qquad= F_{a - b} ^{k + 1} L_n  - F_{m - b} ^{k + 1} L_{n - (m - a)(k + 1)}\,,
\end{split}
\end{equation}
\begin{equation}
\begin{split}
&F_{m - b} \sum\limits_{j = 0}^k {( - 1)^{(a + b + 1)(k - j)} F_{m - a} ^{k - j} F_{a - b} ^j L_{n - (m - b)k - (m - a) + (m - b)j} }\\
&\qquad\qquad\qquad= F_{a - b} ^{k + 1} L_n  - ( - 1)^{(a + b + 1)(k + 1)} F_{m - a} ^{k + 1} L_{n - (m - b)(k + 1)}
\end{split}
\end{equation}
and
\begin{equation}
\begin{split}
&F_{a - b} \sum\limits_{j = 0}^k {( - 1)^{(a + b)j} F_{m - a} ^{k - j} F_{m - b} ^j L_{n - (a - b)k + (m - a) + (a - b)j} }\\
&\qquad\qquad= ( - 1)^{(a + b)k} F_{m - b} ^{k + 1} L_n  + ( - 1)^{a + b + 1} F_{m - a} ^{k + 1} L_{n - (a - b)(k + 1)}\,.
\end{split}
\end{equation}
\begin{theorem}\label{thm.tutbc4x}
The following identities hold for positive integer $k$ and arbitrary integers $a$, $b$, $n$, $m$:
\begin{equation}
\sum\limits_{j = 0}^k {( - 1)^{(a + b + 1)(k - j)} \binom kjF_{m - b}^j F_{m - a}^{k - j} \mathcal{F}_{n - (m - b)k + (a - b)j} }  = F_{a - b}^k \mathcal{F}_n\,,
\end{equation}
\begin{equation}
\sum\limits_{j = 0}^k {( - 1)^{(a + b)j} \binom kjF_{a - b}^j F_{m - a}^{k - j} \mathcal{F}_{n - (a - b)k + (m - b)j} }  = ( - 1)^{a + b} F_{m - b}^k \mathcal{F}_n
\end{equation}
and
\begin{equation}
\sum\limits_{j = 0}^k {( - 1)^j \binom kjF_{a - b}^j F_{m - b}^{k - j} \mathcal{F}_{n + (a - b)k + (m - a)j} }  = ( - 1)^{a + b} F_{m - a}^k \mathcal{F}_n\,.
\end{equation}

\end{theorem}
In particular we have
\begin{equation}
\sum\limits_{j = 0}^k {( - 1)^{(a + b + 1)(k - j)} \binom kjF_{m - b}^j F_{m - a}^{k - j} F_{n - (m - b)k + (a - b)j} }  = F_{a - b}^k F_n\,,
\end{equation}
\begin{equation}
\sum\limits_{j = 0}^k {( - 1)^{(a + b)j} \binom kjF_{a - b}^j F_{m - a}^{k - j} F_{n - (a - b)k + (m - b)j} }  = ( - 1)^{a + b} F_{m - b}^k F_n
\end{equation}
and
\begin{equation}
\sum\limits_{j = 0}^k {( - 1)^j \binom kjF_{a - b}^j F_{m - b}^{k - j} F_{n + (a - b)k + (m - a)j} }  = ( - 1)^{a + b} F_{m - a}^k F_n\,;
\end{equation}
and the corresponding Lucas versions:
\begin{equation}
\sum\limits_{j = 0}^k {( - 1)^{(a + b + 1)(k - j)} \binom kjF_{m - b}^j F_{m - a}^{k - j} L_{n - (m - b)k + (a - b)j} }  = F_{a - b}^k L_n\,,
\end{equation}
\begin{equation}
\sum\limits_{j = 0}^k {( - 1)^{(a + b)j} \binom kjF_{a - b}^j F_{m - a}^{k - j} L_{n - (a - b)k + (m - b)j} }  = ( - 1)^{a + b} F_{m - b}^k L_n
\end{equation}
and
\begin{equation}
\sum\limits_{j = 0}^k {( - 1)^j \binom kjF_{a - b}^j F_{m - b}^{k - j} L_{n + (a - b)k + (m - a)j} }  = ( - 1)^{a + b} F_{m - a}^k L_n\,.
\end{equation}
\subsection{Identities involving Pell, Pell-Lucas and generalized Pell numbers}
\subsubsection{Recurrence relations and related identities}
Since $P_0=0$ and $P_{b-a}=(-1)^{a-b-1}P_{a-b}$ for all $a,b\in\Z$ just like in the Fibonacci case; we find that the Pell relations derived from the identity of Corollary \ref{cor.g869wix}, (with $(G_n)\equiv(P_n)$, the Pell sequence, and $(H_n)\equiv(\mathcal{P}_n)$, the generalized Pell sequence), are identical to those derived in section \ref{sec.fibonacci}. Thus, we have
\begin{equation}\label{eq.w7tubhm}
P_{a-b}\mathcal{P}_{n+m}=P_{m-b}\mathcal{P}_{n+a}-(-1)^{a-b}P_{m-a}\mathcal{P}_{n+b}\,,
\end{equation}
and a couple of special instances:
\begin{equation}\label{eq.gfg8yu3}
P_{n-m}\mathcal{P}_{n+m}=P_n\mathcal{P}_n-(-1)^{n-m}P_m\mathcal{P}_m\,,
\end{equation}
\begin{equation}\label{eq.sfqr94f}
P_{2a}\mathcal{P}_{n+m}=P_{m+a}\mathcal{P}_{n+a}-P_{m-a}\mathcal{P}_{n-a}\,,
\end{equation}
\begin{equation}
P_{2k - 1} \mathcal{P}_{n + m}  = P_{m - 2k} \mathcal{P}_{n + 1}  + P_{m - 1} \mathcal{P}_{n + 2k}\,,
\end{equation}
\begin{equation}\label{eq.taio0c7}
P_{2k} \mathcal{P}_{n + m}  = P_m \mathcal{P}_{n + 2k}  - P_{m - 2k} \mathcal{P}_n
\end{equation}
and
\begin{equation}
P_{2m} \mathcal{P}_{2n}  = P_{n + m} \mathcal{P}_{n + m}  - P_{n - m} \mathcal{P}_{n - m}\,.
\end{equation}
From identity \eqref{eq.gfg8yu3}, we see that Pell numbers also obey Catalan's identity:
\begin{equation}
P_{n-m}P_{n+m}=P_n^2+(-1)^{n+m+1}P_m^2\,.
\end{equation}
We have the following particular cases of identity \eqref{eq.sfqr94f}:
\begin{equation}
P_{2a}P_{n+m}=P_{m+a}P_{n+a}-P_{m-a}P_{n-a}
\end{equation}
and
\begin{equation}
P_{2a}Q_{n+m}=P_{m+a}Q_{n+a}-P_{m-a}Q_{n-a}\,,
\end{equation}
with the special evaluations:
\begin{equation}
2P_{n+m}=P_{m+1}P_{n+1}-P_{m-1}P_{n-1}
\end{equation}
and
\begin{equation}
2Q_{n+m}=P_{m+1}Q_{n+1}-P_{m-1}Q_{n-1}\,.
\end{equation}
\subsubsection{Summation identities}
Summation identities given in Theorems \ref{thm.srgftrf} and \ref{thm.vfhvqcw} are derived by making appropriate substitutions from identity \eqref{eq.w7tubhm} into Lemmata \ref{lem.s9jfs7n} and \ref{lem.binomial}. 
\begin{theorem}\label{thm.srgftrf}
The following identities hold for arbitrary integers $a$, $b$, $n$, $m$ and $k$:
\begin{equation}
\begin{split}
&( - 1)^{a + b + 1} P_{m - a} \sum\limits_{j = 0}^k {P_{m - b} ^{k - j} P_{a - b} ^j \mathcal{P}_{n - (m - a)k - (m - b) + (m - a)j} }\\
&\qquad\qquad\qquad\qquad= P_{a - b} ^{k + 1} \mathcal{P}_n  - P_{m - b} ^{k + 1} \mathcal{P}_{n - (m - a)(k + 1)}\,,
\end{split}
\end{equation}
\begin{equation}
\begin{split}
&P_{m - b} \sum\limits_{j = 0}^k {( - 1)^{(a + b + 1)(k - j)} P_{m - a} ^{k - j} P_{a - b} ^j \mathcal{P}_{n - (m - b)k - (m - a) + (m - b)j} }\\
&\qquad\qquad\qquad= P_{a - b} ^{k + 1} \mathcal{P}_n  - ( - 1)^{(a + b + 1)(k + 1)} P_{m - a} ^{k + 1} \mathcal{P}_{n - (m - b)(k + 1)}
\end{split}
\end{equation}
and
\begin{equation}
\begin{split}
&P_{a - b} \sum\limits_{j = 0}^k {( - 1)^{(a + b)j} P_{m - a} ^{k - j} P_{m - b} ^j \mathcal{P}_{n - (a - b)k + (m - a) + (a - b)j} }\\
&\qquad\qquad= ( - 1)^{(a + b)k} P_{m - b} ^{k + 1} \mathcal{P}_n  + ( - 1)^{a + b + 1} P_{m - a} ^{k + 1} \mathcal{P}_{n - (a - b)(k + 1)}\,.
\end{split}
\end{equation}

\end{theorem}
In particular, we have
\begin{equation}
\begin{split}
&( - 1)^{a + b + 1} P_{m - a} \sum\limits_{j = 0}^k {P_{m - b} ^{k - j} P_{a - b} ^j P_{n - (m - a)k - (m - b) + (m - a)j} }\\
&\qquad\qquad\qquad\qquad= P_{a - b} ^{k + 1} P_n  - P_{m - b} ^{k + 1} P_{n - (m - a)(k + 1)}\,,
\end{split}
\end{equation}
\begin{equation}
\begin{split}
&P_{m - b} \sum\limits_{j = 0}^k {( - 1)^{(a + b + 1)(k - j)} P_{m - a} ^{k - j} P_{a - b} ^j P_{n - (m - b)k - (m - a) + (m - b)j} }\\
&\qquad\qquad\qquad= P_{a - b} ^{k + 1} P_n  - ( - 1)^{(a + b + 1)(k + 1)} P_{m - a} ^{k + 1} P_{n - (m - b)(k + 1)}
\end{split}
\end{equation}
and
\begin{equation}
\begin{split}
&P_{a - b} \sum\limits_{j = 0}^k {( - 1)^{(a + b)j} P_{m - a} ^{k - j} P_{m - b} ^j P_{n - (a - b)k + (m - a) + (a - b)j} }\\
&\qquad\qquad= ( - 1)^{(a + b)k} P_{m - b} ^{k + 1} P_n  + ( - 1)^{a + b + 1} P_{m - a} ^{k + 1} P_{n - (a - b)(k + 1)}\,;
\end{split}
\end{equation}
and the corresponding results involving Pell-Lucas numbers:
\begin{equation}
\begin{split}
&( - 1)^{a + b + 1} P_{m - a} \sum\limits_{j = 0}^k {P_{m - b} ^{k - j} P_{a - b} ^j Q_{n - (m - a)k - (m - b) + (m - a)j} }\\
&\qquad\qquad\qquad\qquad= P_{a - b} ^{k + 1} Q_n  - P_{m - b} ^{k + 1} Q_{n - (m - a)(k + 1)}\,,
\end{split}
\end{equation}
\begin{equation}
\begin{split}
&P_{m - b} \sum\limits_{j = 0}^k {( - 1)^{(a + b + 1)(k - j)} P_{m - a} ^{k - j} P_{a - b} ^j Q_{n - (m - b)k - (m - a) + (m - b)j} }\\
&\qquad\qquad\qquad= P_{a - b} ^{k + 1} Q_n  - ( - 1)^{(a + b + 1)(k + 1)} P_{m - a} ^{k + 1} Q_{n - (m - b)(k + 1)}
\end{split}
\end{equation}
and
\begin{equation}
\begin{split}
&P_{a - b} \sum\limits_{j = 0}^k {( - 1)^{(a + b)j} P_{m - a} ^{k - j} P_{m - b} ^j Q_{n - (a - b)k + (m - a) + (a - b)j} }\\
&\qquad\qquad= ( - 1)^{(a + b)k} P_{m - b} ^{k + 1} Q_n  + ( - 1)^{a + b + 1} P_{m - a} ^{k + 1} Q_{n - (a - b)(k + 1)}\,.
\end{split}
\end{equation}
\begin{theorem}\label{thm.vfhvqcw}
The following identities hold for positive integer $k$ and arbitrary integers $a$, $b$, $n$, $m$:
\begin{equation}
\sum\limits_{j = 0}^k {( - 1)^{(a + b + 1)(k - j)} \binom kjP_{m - b}^j P_{m - a}^{k - j} \mathcal{P}_{n - (m - b)k + (a - b)j} }  = P_{a - b}^k \mathcal{P}_n\,,
\end{equation}
\begin{equation}
\sum\limits_{j = 0}^k {( - 1)^{(a + b)j} \binom kjP_{a - b}^j P_{m - a}^{k - j} \mathcal{P}_{n - (a - b)k + (m - b)j} }  = ( - 1)^{a + b} P_{m - b}^k \mathcal{P}_n
\end{equation}
and
\begin{equation}
\sum\limits_{j = 0}^k {( - 1)^j \binom kjP_{a - b}^j P_{m - b}^{k - j} \mathcal{P}_{n + (a - b)k + (m - a)j} }  = ( - 1)^{a + b} P_{m - a}^k \mathcal{P}_n\,.
\end{equation}

\end{theorem}
In particular we have
\begin{equation}
\sum\limits_{j = 0}^k {( - 1)^{(a + b + 1)(k - j)} \binom kjP_{m - b}^j P_{m - a}^{k - j} P_{n - (m - b)k + (a - b)j} }  = P_{a - b}^k P_n\,,
\end{equation}
\begin{equation}
\sum\limits_{j = 0}^k {( - 1)^{(a + b)j} \binom kjP_{a - b}^j P_{m - a}^{k - j} P_{n - (a - b)k + (m - b)j} }  = ( - 1)^{a + b} P_{m - b}^k P_n
\end{equation}
and
\begin{equation}
\sum\limits_{j = 0}^k {( - 1)^j \binom kjP_{a - b}^j P_{m - b}^{k - j} P_{n + (a - b)k + (m - a)j} }  = ( - 1)^{a + b} P_{m - a}^k P_n\,;
\end{equation}
and the corresponding Pell-Lucas versions:
\begin{equation}
\sum\limits_{j = 0}^k {( - 1)^{(a + b + 1)(k - j)} \binom kjP_{m - b}^j P_{m - a}^{k - j} Q_{n - (m - b)k + (a - b)j} }  = P_{a - b}^k Q_n\,,
\end{equation}
\begin{equation}
\sum\limits_{j = 0}^k {( - 1)^{(a + b)j} \binom kjP_{a - b}^j P_{m - a}^{k - j} Q_{n - (a - b)k + (m - b)j} }  = ( - 1)^{a + b} P_{m - b}^k Q_n
\end{equation}
and
\begin{equation}
\sum\limits_{j = 0}^k {( - 1)^j \binom kjP_{a - b}^j P_{m - b}^{k - j} Q_{n + (a - b)k + (m - a)j} }  = ( - 1)^{a + b} P_{m - a}^k Q_n\,.
\end{equation}
\subsection{Identities involving Jacobsthal, Jacobsthal-Lucas and generalized Jacobsthal numbers}
\subsubsection{Recurrence relations and related identities}
With $(G_n)\equiv(J_n)$, the Jacobsthal sequence, and $(H_n)\equiv(\mathcal{J}_n)$, the generalized Jacobsthal sequence, the identity of Corollary \ref{cor.g869wix} reduces to
\begin{equation}\label{eq.m6yjh0o}
J_{a-b}\mathcal{J}_{n+m}=J_{m-b}\mathcal{J}_{n+a}-(-1)^{a-b}2^{a-b}J_{m-a}\mathcal{J}_{n+b}\,.
\end{equation}
Proceeding as in section \ref{sec.fibonacci}, we have the following particular instances of identity \eqref{eq.m6yjh0o}:
\begin{equation}\label{eq.hk9f7ks}
J_{n-m}\mathcal{J}_{n+m}=J_n\mathcal{J}_n-(-1)^{n-m}2^{n-m}J_m\mathcal{J}_m\,,
\end{equation}
\begin{equation}\label{eq.zcuuakr}
J_{2a}\mathcal{J}_{n+m}=J_{m+a}\mathcal{J}_{n+a}-J_{m-a}\mathcal{J}_{n-a}\,,
\end{equation}
\begin{equation}
J_{2k - 1} \mathcal{J}_{n + m}  = J_{m - 2k} \mathcal{J}_{n + 1}  + J_{m - 1} \mathcal{J}_{n + 2k}\,,
\end{equation}
\begin{equation}\label{eq.qqq235j}
J_{2k} \mathcal{J}_{n + m}  = J_m \mathcal{J}_{n + 2k}  - J_{m - 2k} \mathcal{J}_n
\end{equation}
and
\begin{equation}
J_{2m} \mathcal{J}_{2n}  = J_{n + m} \mathcal{J}_{n + m}  - J_{n - m} \mathcal{J}_{n - m}\,.
\end{equation}
Identity \eqref{eq.hk9f7ks} is a generalization of
\begin{equation}
J_{n-m}J_{n+m}=J_n^2+(-1)^{n+m+1}2^{n-m}J_m^2\,,
\end{equation}
which is the Jacobsthal version of Catalan's identity.
\subsubsection{Summation identities}
Summation identities given in Theorems \ref{thm.nidci3o} and \ref{thm.lvwg18m} are derived by making appropriate substitutions from identity \eqref{eq.m6yjh0o} into Lemmata \ref{lem.s9jfs7n} and \ref{lem.binomial}. 
\begin{theorem}\label{thm.nidci3o}
The following identities hold for arbitrary integers $a$, $b$, $n$, $m$ and $k$:
\begin{equation}
\begin{split}
&( - 1)^{a + b + 1} 2^{a-b}J_{m - a} \sum\limits_{i = 0}^k {J_{m - b} ^{k - i} J_{a - b} ^i \mathcal{J}_{n - (m - a)k - (m - b) + (m - a)i} }\\
&\qquad\qquad\qquad\qquad= J_{a - b} ^{k + 1} \mathcal{J}_n  - J_{m - b} ^{k + 1} \mathcal{J}_{n - (m - a)(k + 1)}\,,
\end{split}
\end{equation}
\begin{equation}
\begin{split}
&J_{m - b} \sum\limits_{i = 0}^k {( - 1)^{(a + b + 1)(k - i)} 2^{(a-b)(k-i)}J_{m - a} ^{k - i} J_{a - b} ^i \mathcal{J}_{n - (m - b)k - (m - a) + (m - b)i} }\\
&\qquad\qquad\qquad= J_{a - b} ^{k + 1} \mathcal{J}_n  - ( - 1)^{(a + b + 1)(k + 1)} 2^{(a-b)(k+1)}J_{m - a} ^{k + 1} \mathcal{J}_{n - (m - b)(k + 1)}
\end{split}
\end{equation}
and
\begin{equation}
\begin{split}
&J_{a - b} \sum\limits_{i = 0}^k {( - 1)^{(a + b)i} 2^{(a-b)(k-i)}J_{m - a} ^{k - i} J_{m - b} ^i \mathcal{J}_{n - (a - b)k + (m - a) + (a - b)i} }\\
&\qquad\qquad= ( - 1)^{(a + b)k} J_{m - b} ^{k + 1} \mathcal{J}_n  + ( - 1)^{a + b + 1} 2^{(a-b)(k+1)}J_{m - a} ^{k + 1} \mathcal{J}_{n - (a - b)(k + 1)}\,.
\end{split}
\end{equation}

\end{theorem}
In particular, we have
\begin{equation}
\begin{split}
&( - 1)^{a + b + 1} 2^{a-b}J_{m - a} \sum\limits_{i = 0}^k {J_{m - b} ^{k - i} J_{a - b} ^i J_{n - (m - a)k - (m - b) + (m - a)i} }\\
&\qquad\qquad\qquad\qquad= J_{a - b} ^{k + 1} J_n  - J_{m - b} ^{k + 1} J_{n - (m - a)(k + 1)}\,,
\end{split}
\end{equation}
\begin{equation}
\begin{split}
&J_{m - b} \sum\limits_{i = 0}^k {( - 1)^{(a + b + 1)(k - i)} 2^{(a-b)(k-i)}J_{m - a} ^{k - i} J_{a - b} ^i J_{n - (m - b)k - (m - a) + (m - b)i} }\\
&\qquad\qquad\qquad= J_{a - b} ^{k + 1} J_n  - ( - 1)^{(a + b + 1)(k + 1)} 2^{(a-b)(k+1)}J_{m - a} ^{k + 1} J_{n - (m - b)(k + 1)}
\end{split}
\end{equation}
and
\begin{equation}
\begin{split}
&J_{a - b} \sum\limits_{i = 0}^k {( - 1)^{(a + b)i} 2^{(a-b)(k-i)}J_{m - a} ^{k - i} J_{m - b} ^i J_{n - (a - b)k + (m - a) + (a - b)i} }\\
&\qquad\qquad= ( - 1)^{(a + b)k} J_{m - b} ^{k + 1} J_n  + ( - 1)^{a + b + 1} 2^{(a-b)(k+1)}J_{m - a} ^{k + 1} J_{n - (a - b)(k + 1)}\,;
\end{split}
\end{equation}
and the corresponding results involving Jacobsthal-Lucas numbers:
\begin{equation}
\begin{split}
&( - 1)^{a + b + 1} 2^{a-b}J_{m - a} \sum\limits_{i = 0}^k {J_{m - b} ^{k - i} J_{a - b} ^i j_{n - (m - a)k - (m - b) + (m - a)i} }\\
&\qquad\qquad\qquad\qquad= J_{a - b} ^{k + 1} j_n  - J_{m - b} ^{k + 1} j_{n - (m - a)(k + 1)}\,,
\end{split}
\end{equation}
\begin{equation}
\begin{split}
&J_{m - b} \sum\limits_{i = 0}^k {( - 1)^{(a + b + 1)(k - i)} 2^{(a-b)(k-i)}J_{m - a} ^{k - i} J_{a - b} ^i j_{n - (m - b)k - (m - a) + (m - b)i} }\\
&\qquad\qquad\qquad= J_{a - b} ^{k + 1} j_n  - ( - 1)^{(a + b + 1)(k + 1)} 2^{(a-b)(k+1)}J_{m - a} ^{k + 1} j_{n - (m - b)(k + 1)}
\end{split}
\end{equation}
and
\begin{equation}
\begin{split}
&J_{a - b} \sum\limits_{i = 0}^k {( - 1)^{(a + b)i} 2^{(a-b)(k-i)}J_{m - a} ^{k - i} J_{m - b} ^i j_{n - (a - b)k + (m - a) + (a - b)i} }\\
&\qquad\qquad= ( - 1)^{(a + b)k} J_{m - b} ^{k + 1} j_n  + ( - 1)^{a + b + 1} 2^{(a-b)(k+1)}J_{m - a} ^{k + 1} j_{n - (a - b)(k + 1)}\,.
\end{split}
\end{equation}
\begin{theorem}\label{thm.lvwg18m}
The following identities hold for positive integer $k$ and arbitrary integers $a$, $b$, $n$, $m$:
\begin{equation}
\sum\limits_{i = 0}^k {( - 1)^{(a + b + 1)(k - i)} \binom kiJ_{m - b}^i 2^{(a-b)(k-i)}J_{m - a}^{k - i} \mathcal{J}_{n - (m - b)k + (a - b)i} }  = J_{a - b}^k \mathcal{J}_n\,,
\end{equation}
\begin{equation}
\sum\limits_{i = 0}^k {( - 1)^{(a + b)i} \binom kiJ_{a - b}^i 2^{(a-b)(k-i)}J_{m - a}^{k - i} \mathcal{J}_{n - (a - b)k + (m - b)i} }  = ( - 1)^{a + b} J_{m - b}^k \mathcal{J}_n
\end{equation}
and
\begin{equation}
\sum\limits_{i = 0}^k {( - 1)^i \binom kiJ_{a - b}^i J_{m - b}^{k - i} \mathcal{J}_{n + (a - b)k + (m - a)i} }  = ( - 1)^{a + b} 2^{(a-b)k}J_{m - a}^k \mathcal{J}_n\,.
\end{equation}

\end{theorem}
In particular we have
\begin{equation}
\sum\limits_{i = 0}^k {( - 1)^{(a + b + 1)(k - i)} \binom kiJ_{m - b}^i 2^{(a-b)(k-i)}J_{m - a}^{k - i} J_{n - (m - b)k + (a - b)i} }  = J_{a - b}^k J_n\,,
\end{equation}
\begin{equation}
\sum\limits_{i = 0}^k {( - 1)^{(a + b)i} \binom kiJ_{a - b}^i 2^{(a-b)(k-i)}J_{m - a}^{k - i} J_{n - (a - b)k + (m - b)i} }  = ( - 1)^{a + b} J_{m - b}^k J_n
\end{equation}
and
\begin{equation}
\sum\limits_{i = 0}^k {( - 1)^i \binom kiJ_{a - b}^i J_{m - b}^{k - i} J_{n + (a - b)k + (m - a)i} }  = ( - 1)^{a + b} 2^{(a-b)k}J_{m - a}^k J_n\,;
\end{equation}
and the corresponding Jacobsthal-Lucas versions:
\begin{equation}
\sum\limits_{i = 0}^k {( - 1)^{(a + b + 1)(k - i)} \binom kiJ_{m - b}^i 2^{(a-b)(k-i)}J_{m - a}^{k - i} j_{n - (m - b)k + (a - b)i} }  = J_{a - b}^k j_n\,,
\end{equation}
\begin{equation}
\sum\limits_{i = 0}^k {( - 1)^{(a + b)i} \binom kiJ_{a - b}^i 2^{(a-b)(k-i)}J_{m - a}^{k - i} j_{n - (a - b)k + (m - b)i} }  = ( - 1)^{a + b} J_{m - b}^k j_n
\end{equation}
and
\begin{equation}
\sum\limits_{i = 0}^k {( - 1)^i \binom kiJ_{a - b}^i J_{m - b}^{k - i} j_{n + (a - b)k + (m - a)i} }  = ( - 1)^{a + b} 2^{(a-b)k}J_{m - a}^k j_n\,.
\end{equation}

\hrule

\noindent 2010 {\it Mathematics Subject Classification}:
Primary 11B39; Secondary 11B37.

\noindent \emph{Keywords: }
Horadam sequence, Fibonacci number, Lucas number, Pell number, Pell-Lucas number, Jacobsthal number, Jacobsthal Lucas number.

\hrule




\begin{thebibliography}{99}
\bibitem{adegoke18} Kunle Adegoke, Weighted sums of some second-order sequences, \emph{The Fibonacci Quarterly} {\bf 56}:3 (2018), 252--262.

\bibitem{aydin} F. T. Aydin,  On generalizations of the Jacobsthal sequence, \emph{Notes on number theory and discrete mathematics} {\bf 24}:1 (2018), 120--135.

\bibitem{halton65} J. H. Halton, On a general Fibonacci identity, \emph{The Fibonacci Quarterly} {\bf 3}:1 (1965), 31--43.

\bibitem{horadam71} A.~F.~Horadam, Pell identities, \emph{The Fibonacci Quarterly} {\bf 9}:2 (1971), 245--252.

\bibitem{horadam96} A.~F.~Horadam, Jacobsthal representation numbers, \emph{The Fibonacci Quarterly} {\bf 34}:1 (1996), 40--54.

\bibitem{koshy} T.~Koshy, \emph{Fibonacci and Lucas numbers with applications}, Wiley-Interscience, (2001).

\bibitem{koshy14} T.~Koshy, \emph{Pell and Pell-Lucas numbers with applications}, Springer Berlin, (2014).

\bibitem{patel13} N.~Patel and P.~Shrivastava, Pell and Pell-Lucas identities, \emph{Global journal of mathematical sciences: theory and practical} {\bf 5}:4 (2013), 229--236.

\bibitem{vajda} S.~Vajda, \emph{Fibonacci and Lucas numbers, and the golden section: theory and applications}, Dover Press, (2008).


\end{thebibliography}
\end{document}